\documentclass{elsart}

\usepackage{amssymb}
\newcommand\toe{\stackrel{e}{\to}}

\begin{document}

\begin{frontmatter}

\title{UPPER BOUND ON THE EDGE FOLKMAN NUMBER $F_e(3,3,3;13)$ }


\author{Nikolay Rangelov Kolev }

\corauth{Supported by the Scientific Research Fund of the St. Kl.
Ohridski Sofia University under contract 226/2010}

\address{Department of Algebra,\\
  Faculty of Mathematics and Informatics,\\
  ``St. Kl. Ohridski'' University of Sofia,\\
  5 J. Bourchier blvd, 1164 Sofia,\\
  BULGARIA}

\begin{abstract}
For a given graph $G$ let $V(G)$ and $E(G)$ denote the vertex and
  the edge set of $G$ respectively.  The symbol $G \toe (a_1, \ldots ,
  a_r)$ means that in every $r$-coloring of $E(G)$ there exists a
  monochromatic $a_i$-clique of color $i$ for some $i \in
   \{1,...,r\}$.  The edge Folkman numbers are defined by the equality
 $$ F_e(a_1, \ldots , a_r ;q) = \min\{ | V(G)| : G \toe  (a_1, \ldots ,
    a_r; q) \mbox{ and } cl(G)<q  \}.
 $$
 It is clear from the definition of edge Folkman numbers  that they
 are a generalization of the classical Ramsey numbers. The problem of
 computation of edge Folkman numbers is extremely difficult and so
 far only eleven edge Folkman numbers  are known.
  In this paper we  prove  the following upper bound  on the number
  $F_e(3,3,3;13)$, namely   $ F_e(3,3,3;13)\leq 30.$  So far it was only known that
 $ F_e(3,3,3;13) < \infty .$
\end{abstract}


\end{frontmatter}

\section{Introduction}

We consider only finite, non-oriented graphs without loops and
multiple edges. We call a $p$-clique of the graph G a set of $p$
vertices, each two of which are adjacent.  The largest positive
integer $p$, such that the graph $G$ contains a $p$-clique is called
a clique number of  $G$ and is denoted by $cl(G)$. We denote by
$V(G)$ and $E(G)$ the vertex set and the edge set of the graph $G $
respectively.  We shall also use the following notations
\begin{itemize}
    \item $N(v)$, $v\in V(G)$ is the set of all vertices of $G$ adjacent to $v$;
  \item $G[V]$, $V \subseteq V(G)$ is the subgraph of $G$ induced by
  $V$;
  \item $G(v) = G[N(v)]$, $v\in V(G)$, that is the subgraph induced by
  the vertices adjacent to $v$  in G;
  \item $K_n$ is the complete graph on $n$ vertices;
   \item $C_n$ is the cycle on  $n$ vertices;
   \item  if $ U_1, U_2 \subseteq V(G) $ then $E(U_1, U_2)$ is the
   set of all edges in $E(G)$ connecting a vertex $ U_1 $ with a vertex of
   $U_2$.
  \end{itemize}
        The Zykov sum of two graphs
        $$G_1+G_2$$
        is the graph obtained
        from the graphs $G_1$ and $G_2$ when we connect each vertex
        from $G_1$ with each vertex from $G_2$.

\textbf{ Definition 1.} Let $a_1, \ldots , a_r$ be positive
integers. The symbol $G \toe (a_1, \ldots ,
a_r)$ means that for each  coloring of the edges of
$G$ in $r$  colors ($r$-coloring) there is a monochromatic  $a_i$-clique in the $i$-th color
for some $i \in \{ 1, \dots , r \}$.

    The Ramsey number $R(a_1,\ldots,  a_r)$  is defined as the least n for
    which $K_n \toe (a_1,\ldots,  a_r)$.

The edge Folkman numbers are defined by the equality
$$ F_e(a_1,\ldots,  a_r ;q) = \min\{ | V(G)| :
  G \toe  (a_1,\ldots,  a_r; q) \mbox{ and } cl(G)<q  \}. $$

It is known that  $    F_e(a_1,\ldots,  a_r;q)$
     exists
    if and only if
  $q> \mbox{max} \{ a_1,\ldots,  a_r \}.$
 This was proved for two colors by Folkman in \cite{Folkman} and in
the general case by Nesertil and Rodl in \cite{NR}.
 It follows from the definition of
$R(a_1,\ldots, a_r)$ that $F_e(a_1,\ldots,  a_r;q)=R(a_1,\ldots,
a_r)$   if  $q>R(a_1,\ldots, a_r).$ In particular $F_e(3,3;q) =6 $
when $q \geq 6 $ because $R(3,3)=6.$ In  1967  P. Erdos       posed
the problems to
 compute  $F_e(3,3;q)$ when $q<6$.
   In  \cite{Graham}   Graham computed the number $F_e(3,3;6)=8$.  He
  established the upper bound proving that $K_3+C_5 \toe (3,3)$.
    An example of a graph G on 15 vertices
  with the properties $G\toe (3,3)$ and   $ cl(G) <5$
  was constructed by Nenov in \cite{N3315}  thus proving that
  $F_e(3,3;5) \leq 15.$
  In \cite{PRU} Piwakowski, Radziszowski,
 Urbanski proved the opposite inequality  $F_e(3,3;5) \geq 15.$
Thus it was proved that $F_e(3,3;5) = 15.$  The last of these
Erdos`s  problems:  to compute the number
  $F_e(3,3;4) $ is still open.
           In \cite{DR}  Dudek and Rodl proved that
          $F_e(3,3;4) \leq 941$.
         The latest  lower bound is $F_e(3,3;4) \geq 19 $
         established by
        S. Radziszowski and Xu Xiaodong
        in \cite{rxiao}.
All these three Erdos problems were about edge Folkman numbers
$F_e(3,3;q)$ which are not equal to
 the Ramsey number $R(3,3)=6.$




Here we shall discuss the edge Folkman numbers  $F_e(3,3,3;q)$ which
are not equal to the Ramsey number $R(3,3,3)=17$. In the most
restricted case we know only the general fact that  $F_e(3,3,3;4)<
\infty.$

 Recently Dudek, Frankl, Rodl \cite{DRF} posed the following

\textbf{ Problem}  Is it true
 that
  $F_e(3,3,3;4)\leq 3^{81}?$

Nenov as consequence of a more general result proved in \cite{Nnew}
that
  $F_e(3,3,3;4) \geq 40.$
 The only
edge Folkman numbers that are not equal to  $R(3,3,3) =17 $  which we
know are:
 $F_e(3,3,3;17)=19 ,$ \cite{Lin} ;
 $F_e(3,3,3;16)=21  $
  (lower bound in \cite{Lin} and upper bound in  \cite{N33316});
 $F_e(3,3,3;15)=23 ,$    \cite{N33315};
 $F_e(3,3,3;14)= 25 ,$  \cite{N33314}.

Define the graph $H = C_5 + C_5 + C_5 + C_5 + C_5 + C_5,$ that is $H$ is a Zykov sum of six copies of the 5-cycle $C_5.$
The main goal of this paper  is to prove the following results:

    \textbf{Theorem} $  H  \toe (3,3,3).$

    As $cl(H) =12$ we obtain the following corollary  from the
    theorem

  \textbf{Corollary}  $F_e(3,3,3;13) \leq 30 .$

    So far it was only known  $   F_e(3,3,3;13) < \infty ,$
which follows from the already cited    general result by Nesertil
and Rodl in \cite{NR} that guarantees the existence of edge Folkman
numbers. The latest lower bound $ 27 \leq F_e(3,3,3;13) $ was
obtained by Nenov  in \cite{Nnew}.

\section{Preliminary results}

Except the graph $H$ that we defined before the theorem we shall
also need the following graphs:
$$S= C_5+ C_5+ C_5+ C_5+ C_5 = H - C_5$$
$$T = K_4 + C_5+ C_5+ C_5+ C_5$$
$$L = \hat{K}_4 + C_5+ C_5+ C_5+ C_5, $$
  where  $\hat{K}_4$ denotes the graph  $K_4$ with one edge deleted.

We shall use the following statement  from \cite{Ndiss}.

 \textbf{Lemma 1} Consider a given disjunct partition  of  $V(T) = V_1 \cup V_2 \cup
 V_3$,
 $ V_i \cap  V_j = \emptyset ,$  $i \neq j,$
such that  $ V_i \cap K_4 \neq \emptyset ,$
 for each  $ i = 1,2,3$.
Then for some  $ i$  we have $ T [V_i] \toe (3,3)$.

Consider  a coloring of the edges of an arbitrary graph  $G$ in
three colors ($3$-coloring). We shall call the colors first, second
and third. For each vertex $v \in V(G)$ we denote by  $N_1(v)$,
$N_2(v)$, $N_3(v)$ its neighbors in first, second and third color
respectively. We shall denote  $ G_i(v) =G[ N_i (v)]$ for $ i =
1,2,3$ and
 $G(v) =G[ N (v)]$.
Now we shall prove the following lemmas.

   \textbf{  Lemma 2 } Consider a 3-coloring of   the edges of an arbitrary graph $G$.

        (a)  $\mbox{   }$  If for some  $v \in V(G)$ and  for some $i=1,2,3$ we have  that $G_i(v) \toe(3,3)$, then   there is a monochromatic triangle in this 3-coloring.

            (b) $\mbox{   }$    If for some  $v \in V(G)$ and  for some $i=1,2,3$ we have that $ cl(G_i) \geq 6$, then   there is a monochromatic triangle in this 3-coloring.

    \textbf{Proof.} (a) Let for example $G_1(v) \toe(3,3)$. If some edge in $G_1$
 is in first color, then  this edge together with the vertex $v$ forms a monochromatic triangle
 in first color. Therefore all edges in $G_1(v)$ are colored in two
 colors only (second and third) and it follows from
 $G_1(v) \toe(3,3)$ that there is a monochromatic   triangle.

                (b) The statement of (b) follows directly from from (a) and the fact $K_6 \toe (3,3)$.

 \textbf{Lemma 3} Consider the graph $Q = K_1 +L= K_1 + \hat{K_4} + C_5 + C_5 + C_5 + C_5.$
   We denote by $w$ the only vertex in
 $K_1$ and by $a$ and $b$   the only
  non-adjacent vertices in
   $\hat{K}_4$. Consider a 3-coloring of the edges of the graph
 $Q,$ such that $E(w, V(\hat{K}_4)) $ contains edges in all the three colors
   and  the edges
  $wa$ and $wb$  are in different colors.  Then there is a monochromatic triangle
  in  this 3-coloring.

 \textbf{Proof} The coloring of the edges of $K_1 +L$ into  three colors
 induces in a natural way a disjunct partition of the vertices of
  $L$ into three sets $V_1 $, $V_2 $, $V_3 $, namely: if the edge  $wx$ is in color $i$,
  then the vertex
  $x$ is in $V_i$, $i=1,2,3.$ We add the edge  $ab$. This completes the graph
   $L $ to the graph $T$.
  Then we have from
 \textbf{ Lemma 1} that
 $T[V_i] \toe (3,3),$ for some $i=1,2,3.$
 As the edges  $wa$ and $wb$ are in different colors then the vertices
  $a$ and $b$ are in different sets  $V_i$.
Thus $T[V_i] = L[V_i] = Q_i(w)$. So $Q_i(w) \toe (3,3)$ and \textbf{Lemma 3} follows from
 \textbf{Lemma 2(a).}

\textbf{Lemma 4} Let $v \in V(H)$ and $S$ be the subgraph of $H$,
induced by the five 5-cycles of $H$ not containing $v$. Assume that there is a  $3$-coloring of the edges of
  $H$ without monochromatic triangles.   Then for every such coloring and
   for each color  $ i$
 we have  $N_i(v) \cap V(S) \neq \emptyset$.

\textbf{Proof} Assume the opposite. Let for  example $N_1(v) \cap
V(S) = \emptyset$. Then for each of the five 5-cycles $C_5$ in $S$
there is an edge of the graph $H$  either  in $N_2(v)\cap V(C_5)   $
or  in $ N_3(v) \cap V(C_5) $.Thus either $N_2(v) $  or $ N_3(v)$
contains $K_6$, which contradicts \textbf{Lemma 2 (b)}.

\section{Proof of the theorem}

Assume the opposite. Consider a coloring of the edges of $H$ in three colors
without a monochromatic triangle.

We shall denote the 5-cycles in the graph  $H$ by $C^{(1)}_5$,
$C^{(2)}_5$, $C^{(3)}_5$, $C^{(4)}_5$, $C^{(5)}_5$, $C^{(6)}_5$.
 We shall first prove the following claims.

\textbf{Claim 1} Each    $C^{(i)}_5 $  is a monochromatic subgraph
of $H$ in the considered coloring.

\textbf{Proof.} Assume the opposite and let  for example
 $C^{(1)}_5 = v_1, v_2, v_3, v_4, v_5, v_1$ is
not a monochromatic subgraph, and the edge  $v_1 v_2$ is in first color,  and the edge
$v_1 v_5$ in  second color. By \textbf{Lemma 4} we have that there
is a vertex  $u_1$ belonging to some of the other five 5-cycles,
such that the edge  $ v_1 u_1$ is in third color. Without loss of
generality we may assume that $u_1 \in C^{(2)}_5.$ Let $u_2$ be a
neighbor of $u_1$ in  $C^{(2)}_5.$ We apply \textbf{Lemma 3} for $K_1 = \{v_1\}$
  and the  subgraph  $L$ induced by the vertices
  $  v_2, v_5 , u_1, u_2 $ and the
5-cycles $C^{(3)}_5$, $C^{(4)}_5$, $C^{(5)}_5$, $C^{(6)}_5$
 (the conditions of \textbf{Lemma 3} are fulfilled
because the edges $v_1 v_2$, $v_1 v_5$, $ v_1 u_1$ are in three
different colors and the vertices $v_2$ and $v_5$ are not adjacent).
According to \textbf{Lemma 3} there is a monochromatic triangle, which is a contradiction.

\textbf{Claim 2} Let  $ v \in C^{(i)}_5$. If  $i \neq j$  then
$E(v,V(C^{(j)}_5))$
cannot contain edges in the both colors different from the color of
$C^{(i)}_5.$

\textbf{Proof.} Assume the opposite and let $v_1  \in V(C^{(1)}_5),$ $C^{(1)}_5$
is monochromatic in first color and
  $ E(v_1, V(C^{(2)}_5))$ contains edges
   $v_1 a $ and  $v_1 b$ which are  in second and third color respectively.
Let  $C^{(1)}_5 = v_1, v_2, v_3, v_4, v_5, v_1$ and
$ C^{(2)}_5 = u_1, u_2, u_3, u_4, u_5,u_1. $ We consider
two cases.

        \textbf{First case.} The vertices  $a$ and $b$ are not adjacent.
  Assume that  $a = u_1$ and $b =u_3.$
   Now the edge
     $v_1 u_1$ is in second color and the edge $v_1 u_3$ is in third
     color. We apply
\textbf{Lemma 3} for   $K_1 =\{ v_1 \}$ and the subgraph $L$
induced by the vertices   $v_2, u_1, u_2, u_3 $ and the 5-cycles
$C^{(3)}_5$, $C^{(4)}_5$, $C^{(5)}_5,$ $C^{(6)}_5$. According to
\textbf{Lemma 3} there is a monochromatic triangle, which is a contradiction.

     \textbf{Second case.}  The vertices $a$ and $b$ are adjacent.
     Let for example $ a=u_1$ and $b=u_2$. Now the edge
     $v_1 u_1$ is in second color and the edge $v_1 u_2$ is in third
     color.
        We shall prove that

         \textit{ $v_1 u_3$ is in the same color in
which is the edge $v_1 u_1$, i.e. in second.}

   Indeed, if we  assume that the edge $v_1 u_3$ is in third color,  then we are in the
    situation of the first case for the vertices
  $u_1$ and $u_3$.
   If we assume that the edge $v_1 u_3$ is in first color, then
 we apply \textbf{Lemma 3} for  $K_1 = \{ v_1 \}$ and the subgraph   $L$
induced by the vertices  $v_2$, $u_1$, $u_2$,$u_3$, and the 5-cycles
$C^{(3)}_5$, $C^{(4)}_5$, $C^{(5)}_5,$ $C^{(6)}_5$.
  It follows from \textbf{Lemma 3} that there is a monochromatic triangle which is a contradiction.
  Thus we proved that $v_1 u_3$ is in second color.

   Analogously we prove that  $v_1 u_5$ is in the same color as the
   edge  $v_1 u_2$, i.e. in third color. Now we apply the first
   case for the vertices $a=u_3$ and $b=u_5$ and thus \textbf{Claim 2} is proved.

\textbf{Claim 3} $\mbox{  }$  If $ C^{(i)}_5$ and $ C^{(j)}_5 $
are in two different colors, then the edges in $E(V(C^{(i)}_5),
V(C^{(j)}_5))$ are in the color different from the colors of $
C^{(i)}_5$ and $ C^{(j)}_5 $.

\textbf{Proof. }  Let for example  $C^{(1)}_5$ is in   first color
and $C^{(2)}_5$ is in second.  Let as above $C^{(1)}_5 = v_1, v_2, v_3, v_4, v_5, v_1$ and
$ C^{(2)}_5 = u_1, u_2, u_3, u_4, u_5,u_1.$
Assume the opposite, i.e. $E(V(C^{(1)}_5)), V(C^{(2)}_5))$  contains at least one edge in first or in second color.
Without loss of generality we may consider that $E(V(C^{(1)}_5)), V(C^{(2)}_5))$ contains an edge in second color
and that this edge is $v_1 u_1$.
  Then it follows from
\textbf{Claim 2} that  $E(v_1,V(C^{(2)}_5))$ contains edges in the
first and second color only. It is not possible
 $E(v_1,V(C^{(2)}_5))$ to contain three edges in second color (otherwise
 the vertex  $v_1$ and an edge from   $C^{(2)}_5$
 form a monochromatic triangle in second color). Hence $E(v_1, V(C^{(2)}_5))$
contains at least three edges in first color. Now we consider the vertex $v_2 $.
  Then according to \textbf{Claim 2} two
cases are possible.

        \textbf{First case.} $E(v_2, V(C^{(2)}_5))$ does not contain edges in third
    color. Now     $E(v_2, V(C^{(2)}_5))$ cannot contain three edges in second color,
     because $C^{(2)}_5$ is in second color and $v_2$ together with an edge in $C^{(2)}_5$
     would form a monochromatic triangle. Therefore
     $E(v_2, V(C^{(2)}_5))$      contains at least three edges in first
     color.
  But    $E(v_1, V(C^{(2)}_5))$ contains at least three edges in first color. Therefore one of the vertices of
    $C^{(2)}_5$ and the edge $v_1 v_2 $ form
  a monochromatic triangle in first color - a
contradiction.

    \textbf{Second case.}  $E(v_2,V(C^{(2)}_5))$ contains at least one edge in third color.
In this situation,  according to \textbf{Claim 2},  $E(v_2,V(C^{(2)}_5))$ does not contain edges in second
    color.
    Therefore  $E(v_2,V(C^{(2)}_5))$ contains either three edges in first color or three edges in third color.
    If
     $E(v_2,V(C^{(2)}_5))$ contains at least three edges in
    first color then having in mind that
   $E(v_1, V(C^{(2)}_5))$ contains at least three edges in first color then one of the vertices of
    $C^{(2)}_5$ and the edge $v_1 v_2$ form a
 monochromatic triangle in first color-a
contradiction. If $E(v_2,
V(C^{(2)}_5))$  contains at least three edges in third color,
as we proved that   $E(v_1, V(C^{(2)}_5))$  contains at least three
edges in first color, then there is  a vertex $u$  in $C^{(2)}_5$,
such that the edge  $v_1 u$
   is in first color, and the edge  $v_2 u$ is in third color.
   Now we apply \textbf{Claim 2} for the vertex
       $u$ and the cycle $C^{(2)}_5$ and we obtain a contradiction.
Now  \textbf{Claim 3} is proved.

According to \textbf{Claim 1}  there are three possible situations:

\textbf{First case.}  There are three 5-cycles  $C^{(i)}_5$  of $H$ that
are monochromatic in three different colors.  Let for example  $C^{(1)}_5$ is in first color, $C^{(2)}_5$
is in second color and  $C^{(3)}_5$ is   in third
color. Without loss of generality we may assume that $ C^{(4)}_5  $ is in third color.
It follows from  \textbf{Claim 3}  that the edges in
$E(V(C^{(1)}_5),V(C^{(3)}_5))$ and $E(V(C^{(1)}_5),V(C^{(4)}_5))$
are in second color, and the edges of $E(V(C^{(2)}_5),V(C^{(3)}_5))$ and
$E(V(C^{(2)}_5),V(C^{(4)}_5))$  are in first color.
As there are no monochromatic triangles in first and second color, then
$E(V(C^{(3)}_5),V(C^{(4)}_5))$ contains edges in third color only.
 Then
any two adjacent vertices in  $C^{(3)}_5$ and any two adjacent
vertices in    $C^{(4)}_5$ induce even a monochromatic 4-clique in
third color, which is a contradiction.

     \textbf{Second  case.}  The  5-cycles of  $H$  are monochromatic in exactly   two different colors.
 Then at least three of the 5-cycles are in one and the same color.
 Let for example   $C^{(1)}_5$ is in first color and
         $C^{(2)}_5$,  $C^{(3)}_5$,  $C^{(4)}_5$ are in second color.
        Then it follows from  \textbf{Claim 3} that the edges in
        $E(V(C^{(1)}_5),V(C^{(2)}_5))$, $E(V(C^{(1)}_5),V(C^{(3)}_5))$ $E(V(C^{(1)}_5),V(C^{(4)}_5))$
are in first color. If $v_1 \in C^{(1)}_5$ then
 $C^{(2)}_5+C^{(3)}_5+C^{(4)}_5$ is contained in $H_3(v_1)$. Thus
  $ K_6 \subseteq H_3(v_1) $, which   contradicts \textbf{Lemma 2
  (b)}.

    \textbf{Third   case.} All the 5-cycles of  $H$ are monochromatic in one
 and the same color, for example first. Let  $v_1 \in C^{(1)}_5$.
    Then it follows from \textbf{Claim 2} that the edges in  $E(v_1,V(C^{(j)}_5))$,
    $j= 2, \ldots ,6$ are at most in two colors, one of which is
    first. As the edges of $C^{(j)}_5$,    $j= 2, \ldots ,6$ are in
    first color then it is impossible $E(v_1,V(C^{(j)}_5))$ to
    contain three edges in first color (otherwise an edge from  $C^{(j)}_5$ and the
     vertex $v_1 $ would form a monochromatic triangle in first color).
      Then it follows from \textbf{Claim 2} that  $E(v_1,V(C^{(j)}_5))$,
    $j= 2, \ldots ,6$ contains at least three edges in second color or at least three edges in third color.
    Then there are at least three 5-cycles among
 $C^{(j)}_5$,    $j= 2, \ldots ,6$, such that $E(v_1,V(C^{(j)}_5))$
contains either three edges in second color or three edges in third
color. Thus at least three of the sets  $N_2(v_1) \cap V(C^{(j)}_5)$ contain an edge
or at least three of the sets $N_3(v_1) \cap V(C^{(j)}_5)$ contain an edge.
Therefore  $cl(H_3(v_1) )\geq 6$ or $cl(H_2(v_1))\geq 6$,
which contradicts \textbf{Lemma 2 (b)}. This completes the proof of
the theorem.

    \textbf{Acknowledgement} I am indebted to prof. Nenov whose
    remarks improved the presentation of the paper.

\end{document}